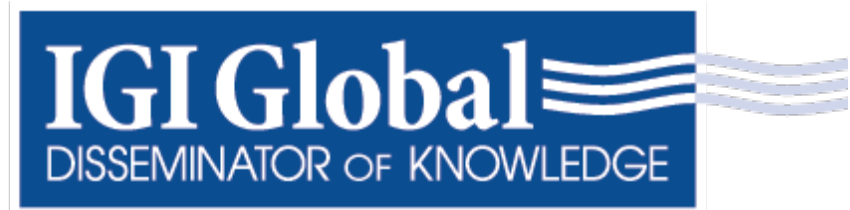

# Testing exchangeability with martingale for change-point detection

Liang Dai[(1)],  Mohamed-Rafik Bouguelia[(2)]
   (1) Ekkono Solutions, Borås, Sweden.
   (2) Halmstad University, Sweden – Department of Intelligent Systems and Digital Design.

**ABSTRACT**

*This work proposes a new exchangeability test for a random sequence through a martingale based approach. Its main contributions include: 1) an additive martingale which is more amenable for designing exchangeability tests by exploiting the Hoeffding-Azuma lemma; 2) different betting functions for constructing the additive martingale are studied. By choosing the underlying probability density function of p-values as a betting function, it can be shown that, when a change-point appears, a satisfying trade-off between the smoothness and expected one-step increment of the martingale sequence can be obtained. An online algorithm based on Beta distribution parametrization for constructing this betting function is discussed in detail as well.*

Keywords: Change Detection, Data Streams, Martingale, Exchangeability Test, Online Learning

## 1. INTRODUCTION

Today, with the big data phenomenon, considerable research is being put on designing online algorithms capable of detecting changes in the distribution of streaming data. Such changes happen in any system due to unforeseen external conditions and can indicate possible faults, anomalies, or any significant events that require the attention of a human expert. A change happens when new data observations deviate from previous ones and become not identically distributed (i.e. way say that the exchangeability assumption is violated). To detect such changes, it is often necessary to model the system's behavior (e.g. have knowledge about the distribution of the data). However, in a time where the volumes of data and the complexity of systems is continuously growing, it becomes not always feasible for human experts to build an exhaustive or good definition of each system's behavior. This paper discusses the exchangeability test for a given random sequence using the martingale-based approach. The results can be applied for detecting change-point in streaming data including cases for abrupt change and concept drift, i.e. changes that happen in a smooth and incremental manner. In real world applications, systems often exhibit complex behavior and the data distribution is unknown. The martingale based approach is well suited in these scenarios for change-point detection since it does not rely on any distributional knowledge about the data, which is in contrast to the traditional sequential change-point detection methods such as Sequential

Probability Ratio Test **[Wald 1945]**, and the Cumulative SUM control chart **[Page 1954]**. The work in this paper can be beneficial to a wide range of applications, such as, monitoring systems for fault detection **[Farouq 2020; Smolyakov 2018]**, detection changes in driving behavior **[Ali 2017]**, emotion change from speed **[Zhaocheng 2016],** and detecting anomalous situations in ambient assisted living applications **[Bersch 2013]**.

The idea of using martingale for exchangeability test dates back to the work in **[Vovk 2003]**, building on the theory of Transductive Confidence Machine **[Vovk 2005]**, where the concept of "exchangeability martingales" was introduced for implementing the test which works in an online manner. Later, it was further established that **[Fedorova 2012]** to maximize the logarithmic growth rate of the multiplicative martingale, the betting function used for constructing the martingale should be chosen as the empirical probability density function (p.d.f.) of the p-values. The way to construct this empirical p.d.f. suggested therein was to use a modified kernel density estimator. In **[Ho 2005a; Ho 2005b]**, the authors applied various concentration inequalities on the multiplicative martingale sequence to design tests for detecting change in data streams. However, due to the high variability of the multiplicative martingale sequence, it is hard to design the test based on concentration inequalities. In addition, the multiplicative martingale values in the log scale will exhibit an undesirable behavior that a decreasing trend even when no change happens (which will be illustrated by an example in the Evaluation section).

In order to address these shortcomings as stated before, this paper proposes a new type of martingale for the exchangeability test (called "*additive martingale*") which can be implemented in an online fashion. Different betting functions for constructing the additive martingale are discussed as well. Interestingly, similar to the multiplicative martingale case, it is shown that by choosing the betting function as the underlying p.d.f. of the p-values, when a change-point appears (then the generated p-values are not uniformly distributed), a satisfying balance between the smoothness and expected one-step increment in the martingale sequence will be obtained. Based on Beta distribution parametrization, a computationally efficient way for constructing this betting function is discussed. Finally, designing tests for change-point detection based on different concentration inequalities are discussed. The novelty and contributions of this paper are summarized as follows:
- Formulation of a new martingale for change-point detection, which makes it possible to define an upper bound confidence using the Hoeffding-Azuma lemma.
- A proof that constructing the martingale sequence using the "underlying p.d.f of the p-values" as a betting function, leads to a satisfying trade-off between: (i) "*how fast the sequence values grow when a change in the data distribution occurs*", and (ii) "*the smoothness of the sequence when no change occurs*".
- An online/incremental algorithm based on the Beta distribution parametrization for constructing the proposed martingale.

The remaining of this paper is organized as follows. Section 2 provides essential definitions and discusses existing martingale-based approaches. Section 3 presents our proposed martingale framework for change-point detection and discusses various betting functions for constructing it. Section 4 shows how to use the Beta distribution parametrization to construct the proposed martingale sequence in an online fashion. Section 5 discusses how to use the Hoeffding-Azuma lemma with the proposed additive martingale in order to design precise tests for change-point detection. Section 6 showcases properties of the proposed method through experimental evaluation. Finally, a conclusion and future work are provided in section 7.

## 2. BACKGROUND

**Definition 1** (Martingale). A sequence of random variables $\{S_n\}_{n=1}^{\infty}$ is a martingale if for any $n \geq 1$, it satisfies that

$$\boldsymbol{E}(S_{n+1}|S_n, S_{n-1}, \cdots, S_1) = S_n \qquad (1)$$

where $E(.)$ refers to the *expectation*.

**Definition 2** (Exchangeability). A set of random variables $Z_1, Z_2, \cdots, Z_n$ are exchangeable if it holds that
$$P(Z_1, Z_2, \cdots, Z_n) = P(Z_{\sigma(1)}, Z_{\sigma(2)}, \cdots, Z_{\sigma(n)}), \quad (2)$$
in which $\sigma(\cdot)$ denotes any permutation of $[1,2,\cdots,n]$. A series of random variables $(Z_1, Z_2, \cdots)$ is exchangeable if $(Z_1, Z_2, \cdots, Z_n)$ is exchangeable for any natural number n.

Let $(z_1, z_2, \cdots)$ denote a sequence of data samples. For each sample $z_i$, the "nonconformity measure" quantifies the strangeness of $z_i$ with respect to the other data samples:
$$\alpha_i = A(z_i, \{z_1, z_2, \cdots, z_n\}). \quad (3)$$

The operator $A(\cdot,\cdot)$ in eq. (3) represents certain algorithm which takes $z_i$ and the other data samples as inputs and returns a value $\alpha_i$ reflecting the "nonconformity" of $z_i$ with respect to the other data samples. For example, one way to obtain the nonconformity measure is based on the Nearest Neighbor algorithm as follows:
$$\alpha_i = min_{j \neq i, j \in [1:n]} d(z_i, z_j),$$
where $d(\cdot,\cdot)$ denotes the Euclidean distance.

Once the nonconformity measures of all data samples are calculated, the sequence of p-values can be calculated using Algorithm 1 **[Fedorova 2012]**, in which $\theta_i$ denotes a random number uniformly distributed in $[0,1]$ and $|A|$ denotes the cardinality of set A.

**Algorithm 1**: p-value calculation for data samples

```
Input: sequence of data (z₁, z₂, ⋯)
Output: sequence of p-values (p₁, p₂, ⋯)
For i = 1, 2, ⋯ do
    Obtain sample zᵢ
    For j = 1, ⋯, i do
        αⱼ = A(zⱼ, {z₁, ⋯, zᵢ})
    EndFor
```
$$p_i = \frac{|\{j: \alpha_j > \alpha_i\}| + \theta_i |\{j: \alpha_j = \alpha_i\}|}{i}$$
```
EndFor
```

**Algorithm 2**: Inductive version of Algorithm 1.

```
Input: Training set (z₋ₘ, z₋ₘ₊₁, ⋯, z₀), sequence of data (z₁, z₂, ⋯)
Output: sequence of p-values (p₁, p₂, ⋯)
For i = 1, 2, ⋯ do
    Obtain sample zᵢ
    αᵢ = A(zᵢ, {z₋ₘ, z₋ₘ₊₁, ⋯, z₀})
```
$$p_i = \frac{|\{j: \alpha_j > \alpha_i\}| + \theta_i |\{j: \alpha_j = \alpha_i\}|}{i}$$
```
EndFor
```

The following Theorem 1 from **[Vovk 2003; Fedorova 2012]** plays a pivotal role for developing the martingale based test for exchangeability.

**Theorem 1**: If the data samples $\{z_1, z_2, \cdots\}$ satisfy the exchangeability assumption, Algorithm 1 will produce p-values $\{p_1, p_2, \cdots\}$ that are independent and uniformly distributed in [0,1].

In Theorem 1, the $p$ values reflect the strangeness of the corresponding data points - a smaller $p$ means a larger strangeness of the corresponding data sample. Note that computing p-values according to Algorithm 1 is heavily time consuming: whenever a new sample is obtained, the nonconformity measures for all the previous data samples must be recalculated. To avoid these expensive computations, the "inductive" version is applied to compute the p-values, as given in Algorithm 2. The "inductive" version assumes a prefixed training set, and based on this fixed training set, the nonconformity values for all the samples only need to be calculated once **[Vovk 2005; Denis 2017]**, hence it is much more computationally efficient.

To prepare for the next few sections, the main idea of the martingale based approach for exchangeability test is briefly summarized here. When the p-value sequence is obtained by running Algorithm 1, a new sequence can be constructed through a proper "betting function" (which satisfy some special properties). If no change happens in the distribution of the data, the newly constructed sequence will be very likely to stay in a bounded region since it is a valid martingale; otherwise the sequence will have a growing or decreasing trend as the p-values will not be uniformly distributed in [0,1] anymore (as implied by Theorem 1, lack of exchangeability will give non-uniformly distributed p-values), and will start concentrating around a small region.

## 2.1. Multiplicative Martingale

In **[Vovk 2003; Ho 2005b; Fedorova 2012]**, the authors proposed the exchangeable martingale (referred to as multiplicative martingale in this work) for the exchangeability test. The idea is summarized as follows. For the sequence $\{p_1, p_2, \cdots\}$ generated by Algorithm 1, consider the following random sequence

$$S_n = \prod_{i=1}^{n} f_i(p_i), \quad n = 1, 2, \cdots, \quad (4)$$

where $f_i(p): [0,1] \to [0, \infty)$ is called "*betting function*", which satisfies

$$\int_0^1 f_i(p_i) dp_i = 1.$$

From which, it follows that

$$E(S_{n+1} | S_n, S_{n-1}, \cdots, S_1) = S_n \quad (5)$$

Therefore $\{S_n\}_{n=1}^{\infty}$ is a valid martingale sequence according to the definition. Different betting functions have been suggested for multiplicative martingales; three typical ones are summarized as follows.

*Power Martingale:* It uses a fixed power function as betting function $f(p) = \epsilon p^{\epsilon-1}$, where $\epsilon \in [0,1]$. Therefore, the power martingale for a given $\epsilon$ is written as

$$S_n^{(\epsilon)} = \prod_{i=1}^{n} \epsilon p_i^{\epsilon-1}. \quad (6)$$

*Mixture Power Martingale*: It uses a mixture of power martingales based on different $\epsilon \in [0,1]$ values.

$$S_n = \int_0^1 S_n^{(\epsilon)} d\epsilon. \quad (7)$$

*Plug-In Martingale*: It uses an empirical p.d.f. of the p-values as betting function. In addition, it has been justified in **[Fedorova 2012]** that the plug-in martingale is more efficient in terms of rapid change in the martingale value when change-point happens. To construct the empirical p.d.f., a modified kernel density estimator is used therein.

Due to the unboundness of the power function and the multiplicative construction, it is found inconvenient for the multiplicative martingale to adapt the Hoffding-Azuma type concentration inequalities to design tests for detecting change in data streams, see the discussions in **[Ho 2005a; Ho 2005b]**. These shortcomings of the multiplicative martingale motivate us to propose the "additive martingale" as an alternative. This will be elaborated further in the next section.

## 3. ADDITIVE MARTINGALE

As pointed out in **[Denis 2017]** "*it is interesting whether there are any other test exchangeability martingales apart from the conformal martingales (i.e. the one defined in eq. (4))*". To address some of the issues of the multiplicative martingale, a new additive martingale is presented in this section.

### 3.1. Basic Idea

The additive martingale is inherently related to the multiplicative martingale, and their connection can be elucidated through the following reasoning. Suppose that the logarithm operation is taken on both sides of eq. (4):

$$\log(S_n) = \sum_{i=1}^{n} \log(f_i(p_i)), \quad n = 1, 2, \cdots. \quad (8)$$

What we hope to get is that $\{\log(S_n)\}_{n=1}^{\infty}$ will be a valid martingale sequence (since we want to obtain a martingale in the "additive" sense), or equivalently, it satisfies that

$$\int_0^1 \log(f_i(p)) \, dp = 0, \quad i = 1, 2, \cdots. \quad (9)$$

However, in the multiplicative martingale case, the betting function is chosen to satisfy

$$\int_0^1 f_i(p) dp = 1, \quad i = 1, 2, \cdots, \quad (10)$$

and $f_i(p) \geq 0$. Note that since the $\log$ function is concave, we will get

$$\int_0^1 \log(f_i(p)) \, dp \leq \log \int_0^1 f_i(p) dp = 0,$$

which gives that eq. (10) will not be able to imply eq. (9).

***Remark***: *The previous reasoning proves that, instead of being a martingale, the sequence $\{\log(S_n)\}_{n=1}^{\infty}$ is actually a supermartingale. This justifies the decreasing trend of the multiplicative martingale value (in the log scale), as illustrated in Figure 2.d.*

To mitigate this problem (i.e. to get a valid martingale), we can directly enforce the betting functions $f_i(p)$ to integrate to zero, that is, for $n = 1, 2, \cdots$, let $S_n$ be defined as

$$S_n = \sum_{i=1}^{n} f_i(p_i) \quad \text{with} \quad \int_0^1 f_i(p) dp = 0. \quad (11)$$

Then we have:

$$E(S_{n+1}|S_n, \cdots, S_1)$$
$$= \int_0^1 \left( \sum_{i=1}^{n} f_i(p_i) + f_{n+1}(p_{n+1}) \right) dp_{n+1} = \sum_{i=1}^{n} f_i(p_i) + \int_0^1 f_{n+1}(p) dp = \sum_{i=1}^{n} f_i(p_i) = S_n,$$

therefore $S_n$ becomes a valid martingale.

### 3.2. Betting functions for additive martingale

In what follows, two betting function constructions are given to get valid additive martingales.

**Shifted odd functions**
By definition, any odd function $g(p): [-1/2, 1/2] \to \mathbb{R}$ will satisfy

$$\int_{-1/2}^{1/2} g(p)dp = 0,$$

from which, it follows that

$$\int_0^1 g(p - 1/2)dp = 0.$$

This simple fact implies that $f(p) = g(p - 1/2)$ will be a valid betting function for any odd function $g(p)$. One example is to let $g(p) = -p$, more betting functions can easily be constructed by picking different odd functions.

**Shifted empirical probability density function**
From the p-values calculated by Algorithm 1, an empirical probability density function of the p-values can be obtained (one computationally efficient way will be discussed later on), which we denote as $\hat{\rho}_t(p)$ at time $t$. Based on this, it can be readily checked that a valid betting function can be formulated as

$$f_t(p) = \hat{\rho}_t(p) - 1. \qquad (12)$$

This construction is not only valid, and in fact, it will give a rapid and smooth change in the martingale sequence when a change-point happens (see the experiment result in Figure 2.b) in the data sequence. To explain this observation, the following optimization problem is defined:

$$max_{g(p)} \left( \int_0^1 g(p)\rho(p)dp \right)^2 - \lambda \int_0^1 g^2(p)dp$$
$$\text{s.t.} \quad \int_0^1 g(p)dp = 1 \text{ and } g(p) \geq 0. \qquad (13)$$

The objective function in eq. (13) consists of two parts: the first part represents the expected increment of the martingale sequence value at each step, when betting function $g(p) - 1$ is used, and given the underlying p.d.f. of p-values as $\rho(p)$; the second term represents the "flatness" (or "regularness") of the betting function.

To make better sense of the optimization problem, the following two extreme cases are analyzed: when $\lambda = 0$ and $\lambda = \infty$.

When $\lambda = 0$, since $g(p)$ and $\rho(p)$ are both non-negative functions, the problem in eq. (13) can be reduced to

$$max_{g(p)} \int_0^1 g(p)\rho(p)dp$$
$$s.t. \quad \int_0^1 g(p)dp = 1 \text{ and } g(p) \geq 0. \qquad (14)$$

Assume that $\rho(p)$ is upper-bounded by $M$ and at point $p_0$ we have $\rho(p_0) = M$, then it follows that
$$\int_0^1 g(p)\rho(p)dp \leq \int_0^1 M g(p)dp = M, \qquad (15)$$
and the maximum value $M$ can be obtained when $g(p)$ is set to be $\delta(p - p_0)$, where $\delta(.)$ denotes the Dirac delta function which is an extremely peaky. More concretely, when $g(p) = \delta(p - p_0)$, we have
$$\int_0^1 g(p)\rho(p)dp = \int_0^1 \rho(p)\delta(p - p_0)dp = \rho(p_0) = M.$$

Figure 1 illustrates how the martingale $S_t$ changes over time, when the betting function is a Dirac delta function (a Gaussian pdf with a very small variance). As shown on Figure 1.b, the martingale can reach very high values when p-values are not uniformly distributed. However, even when p-values are uniformly distributed, the martingale sequence still have a high variation and may end up far from its initial point, as can be observed from Figure 1.a. This is not ideal for change-point detection, since it may increase the possibilities of false-alarms.

Let's discuss the case when $\lambda = \infty$. In this situation, the problem in eq. (13) reduces to
$$min_{g(p)} \int_0^1 g^2(p)dp$$
$$s.t. \int_0^1 g(p)dp = 1 \text{ and } g(p) \geq 0. \qquad (16)$$

Given by Cauchy-Schwarz inequality, we have that
$$1 = \left(\int_0^1 g(p)dp\right)^2 \leq \int_0^1 1^2 dp \int_0^1 g^2(p)dp = \int_0^1 g^2(p)dp,$$
where the equality holds when $g(p) = 1$, $p \in [0,1]$. Given by these calculations, in the case when $\lambda = \infty$, the optimal solution to eq. (13) is given by a uniform distribution function within the interval $[0,1]$ - which is the most "regular" function.

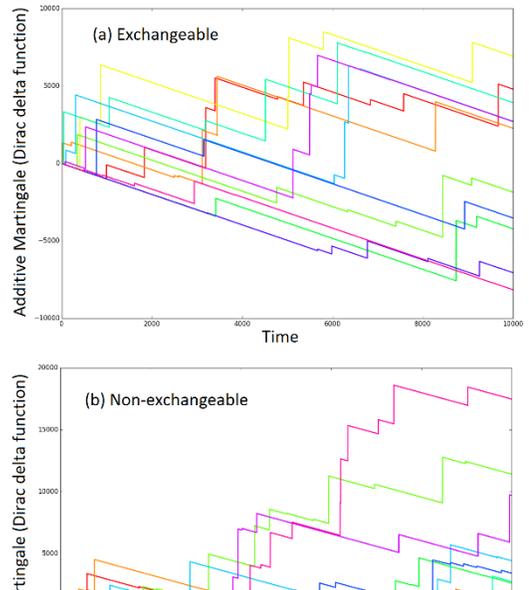

Figure 1: Values of the additive martingale $S_i$ over time when using a simulated Dirac delta function (approximated by Gaussian with very small variance), in two cases: (a) exchangeable and (b) non-exchangeable. Colors represent different runs of the simulation.

The previous discussion implies that a proper choice of $\lambda$ will give a satisfying balance between the one-step increment and the "regularness" (which means that the sequence does not include big jumps) of the martingale sequence. Next, we show that, when choosing
$$\lambda = \int_0^1 \rho^2(p)dp,$$
the optimal solution to eq. (13) is given by $\rho(p)$, which gives that the corresponding betting function is $\rho(p) - 1$.

Note that when $\lambda$ is chosen as $\int_0^1 \rho^2(p)dp$, again by the Cauchy-Schwarz inequality, we will have
$$\left(\int_0^1 g(p)\rho(p)dp\right)^2$$
$$\leq \int_0^1 \rho^2(p)dp \int_0^1 g^2(p)dp, \qquad (17)$$
where the equality holds when $g(p) = c\rho(p)$ and $c$ is a constant. When the equality holds, the maximum of the objective (which is zero) in eq. (13) is achieved. Given the fact that both $g(p)$ and $\rho(p)$ integrate to 1 in [0,1], we get $c = 1$, and it implies that the optimum to the optimization problem in eq. (13) is $\rho(p)$.

## 4. ESTIMATING P-VALUES DISTRIBUTION WITH THE BETA DISTRIBUTION

Given the importance of the p-value density function in constructing efficient additive martingales, in what follows, a computationally efficient way to build up an approximation of the p-value density function is discussed.

According to Theorem 1, when a change-point happens, the p-values will not be uniformly distributed within [0,1]. Typical cases are that the distribution will be skewed with a single mode. This observation inspires us to model the p-value density function with a Beta distribution, which is defined as follows.

### 4.1. Beta distribution

**Definition 3**: The beta distribution $B_{\alpha,\beta}(x)$, parametrized by two positive shape parameters α and β, defines a family of continuous probability distributions on [0,1], given as
$$B_{\alpha,\beta}(x) = \frac{x^{\alpha-1}(1-x)^{\beta-1}}{Beta(\alpha,\beta)},$$
where $Beta(\alpha,\beta) = \frac{\Gamma(\alpha)\Gamma(\beta)}{\Gamma(\alpha+\beta)}$ and $\Gamma(\cdot)$ denotes the Gamma function.

Note that when $\alpha = \beta = 1$, it gives the uniform distribution on [0,1]. When both α and β are greater than one, an imbalanced choice of α, β will give a skewed density function, which is of particular interest to us since it will be useful to model the skewed p-value distribution with a single mode.

***Remark***: *There exist non-parametric approaches **[Tsybakov 2009]** for estimating density functions, for example the histogram and kernel based density estimators. For these estimators, optimal choice for the number of bins or the kernel bandwidth will depend on knowledge of the underlying p.d.f. which is often unknown. The Beta parametric approach presents an alternative for the case when single mode appears in the p-value distribution. The parameters are easy to estimate and the estimation can be done in an online fashion and will be explained in next part.*

### 4.2. Beta distribution

In **[Bain 1992]**, a moment-matching based method for the shape parameters (α and β) estimation was proposed. A notable feature of this approach is that it only involves the calculation of the sample mean and variance.

More concretely, assume we are given a set of p-values as $p_1, \cdots, p_n$, then according to **[Bain 1992]**, $\widehat{\alpha_n}$ and $\widehat{\beta_n}$ (the estimated parameters) can be calculated as follows:
$$\widehat{\alpha_n} = \overline{p_n}\left(\frac{\overline{p_n}(1-\overline{p_n})}{s_n} - 1\right) \quad (18)$$
$$\widehat{\beta_n} = (1-\overline{p_n})\left(\frac{\overline{p_n}(1-\overline{p_n})}{s_n} - 1\right), \quad (19)$$
in which $\overline{p_n}$ and $s_n$ denote the sample mean and variance respectively, which are defined as
$\overline{p_n} = \frac{1}{n}\sum_{i=1}^{n} p_i$ and $s_n = \frac{1}{n-1}\sum_{i=1}^{n}(p_i - \overline{p_n})^2$.

Online calculation of $\overline{p_n}$ and $s_n$ can be found in **[Welford 1962]**. Concretely,
$$\overline{p_n} = \overline{p_{n-1}} + \frac{p_n - \overline{p_{n-1}}}{n}.$$

And for updating the sample variance, one can first update $M_n = \sum_{i=1}^{n}(p_i - \overline{p_n})^2$ recursively as follows
$$M_n = M_{n-1} + (p_n - \overline{p_{n-1}})(p_n - \overline{p_n}), \quad (20)$$
which gives $s_n = \frac{M_n}{n-1}$.

**Remark**: *Note that when sliding window (with size $W$) is introduced, the sample mean $\overline{p_n^W}$ and sample variance $s_n^W$, defined as $\overline{p_n^W} = \frac{1}{W}\sum_{i=n-W+1}^{n} p_i$ and $s_n^W = \frac{1}{W-1} M_n^W$ in which $M_n^W = \sum_{i=n-W+1}^{n}(p_i - \overline{p_n^W})^2$, can similarly be calculated in an online manner, and the main steps to do this can be found in the Appendix.*

## 5. DESIGNING TESTS FOR CHANGE-POINT DETECTION

In this section, the Hoeffding-Azuma inequality and the Doob-Kolmogorov's inequality are applied to the additive martingale sequence to develop statistical tests for change-point detection. The general idea is that, when no change-point appears, the martingale sequence will be bounded in certain region with high probability. However, when the sequence exceeds the specific region, it is very likely that a change-point has occurred, hence an alarm needs to be triggered.

**Theorem 2** (Hoeffding-Azuma inequality):
Let $c_1, \cdots, c_m$ be constants and let $Y_1, \cdots, Y_m$ be a martingale difference sequence with $a_k \leq Y_k \leq b_k$ for each $k$. Then for any $t \geq 0$, we have
$$P\left(\left|\sum_{k=1}^{m} Y_k\right| \geq t\right) \leq 2\exp\left(-\frac{t^2}{\sum_{k=1}^{m}(b_k - a_k)^2}\right) \quad (21)$$

In the additive martingale case, when the betting function is chosen as a shifted odd function $f(p)$ with $|f(p)| \leq 1$, eq. (21) reduces to
$$P(|S_m| \geq t) \leq 2\exp\left(-\frac{t^2}{2m}\right),$$
where $S_m = f(p_1) + \cdots + f(p_m)$.

This fact can be used to design statistical tests for $H_0$ ($H_0$: $p$-values follow a uniform distribution in $[0,1]$, i.e. no change-point appearing). More specifically, given the significance level $\alpha$, when
$$|S_m| > \sqrt{2m \ln\left(\frac{2}{\alpha}\right)},$$
the hypothesis $H_0$ is rejected and an alarm is triggered.
In practice, sliding window (assume the window size is $W$) can be introduced to track the change more rapidly. Given the significance level $\alpha$ and after similar calculations as done before, when
$$|S_m - S_{m-W}| > \sqrt{2W \ln\left(\frac{2}{\alpha}\right)}, \quad (22)$$
the hypothesis $H_0$ will be rejected with the significance level $\alpha$.

**Remark**: *It can be observed that, due to the additive structure of the constructed martingale (in contrast to the multiplicative martingale case), it becomes more convenient to apply the hoeffding-Azuma inequality for the test design.*

Inspired by **[Ho 2005b]**, tests can also be designed based on the following inequality.

**Theorem 3** (Doob-Kolmogorov inequality):
Let $Y_1, Y_2, \cdots, Y_n$ be a martingale difference sequence, and $S_k = Y_1 + \cdots + Y_k$ for $k = 1, \cdots, n$. Then it follows that

$$P\left(\max_{1 \leq k \leq n} |S_k| \geq t\right) \leq \frac{E[S_n^2]}{t^2}. \quad (23)$$

Notice the facts in Theorem 1 that the $p$-values generated from Algorithm 1 are independent from each other, therefore if the betting function $f_i(p)$ is chosen independent from $p_j$, where $j \neq i$, for instance when $f(p) = -p + 1/2$, by calculating out $E[S_n^2]$, the inequality in eq. (23) can be reduced to (assuming a sliding window with size $W$ is used):

$$P\left(\max_{n-W+1 \leq k \leq n} |S_k| \geq t\right) \leq \frac{W}{12t^2},$$

which gives that $H_0$ will be rejected (with significance level $\alpha$) when

$$\max_{n-W+1 \leq k \leq n} |S_k| \geq \sqrt{\frac{W}{12\alpha}}.$$

## 6. EVALUATION

To showcase properties of the proposed method, an experiment is conducted to illustrate how both martingale sequences behave in the exchangeable and non-exchangeable cases. The setup is given as follows: data in the first part (before the vertical dashed line in Figure 2) of the sequence are i.i.d. drawn from a Gaussian distribution with zero mean and unit variance; data in the second part is also drawn from a Gaussian distribution, but with a different mean from the first Gaussian. The second part of the sequence models the data when change has happened. Distance to nearest neighbor is used as the non-conformity measure, and the Algorithm 2 is used for p-value calculation.

The results are reported in Figure 2. In all these figures, the x-axis represents time and the y-axis represents the martingale values (in multiplicative martingale case, it's the martingale values in the log scale). Comparing Figure 2 (a-b), one can observe that the curves obtained with the estimated p.d.f. as the betting function can achieve higher martingale value (which will make the change-point detection more confident) than the one with shifted odd function as the betting function. Again, the reason is that, the approach based on the estimated p.d.f. can adapt to the change and gain a higher one-step increment in the curve. In Figure 2 (b), one can also observe that, when a change in the data distribution just happens, the curve will increase slowly in the initial time steps, during which the algorithm is gaining knowledge of the changed p.d.f. of the p-values, but after this phase, the curves increase very rapidly. In the Figure 2 (c-d) obtained using the multiplicative martingales, the curves will exhibit decreasing trend in the period of no change; after the change appears, the curves start to increase and will take significant time to return to a high value. Though it is possible to post-process the curve, for example by applying a one-step finite difference filter to transform the "v"-shape curve into a step-alike curve similar to the ones in Figure 2 (a-b), or to use another trick introduced in **[Denis 2017]**, however these post-processing will introduce additional complications for designing tests (for example, consider applying the Hoeffding-Azuma type inequalities on the transformed sequence). In addition, in the period of no deviation, the curves in Figure 2 (b) have small variations as compared to the curves in other subfigures, which will make it less prone to trigger false alarms.

The usefulness of the proposed method is also demonstrated in a real-world application scenario related to fault detection. To do so, the proposed method has been applied to real heat-pump data generated over a couple of years, in order to detect potential compressor failures. The data is a multivariate time series consisting of the "ratio of compressor runs", the "additional heat supplied", and the "produced hot water temperature". As shown on Figure 3, the data is subject to change over time due to variations in the outside temperature. The black dots on Figure 3 and the black subsequences on the time series of Figure 4, correspond to a period where a compressor failure occurred. The red curve in Figure 4 shows the values of the martingale sequence (i.e. the deviation level), and the green dots on the same figure show the p-values. One can observe that the p-values became close to 0 during the time where the compressor failure happened while the red curve quickly increased, correctly indicating an abnormal change in the heat-pump's behavior.

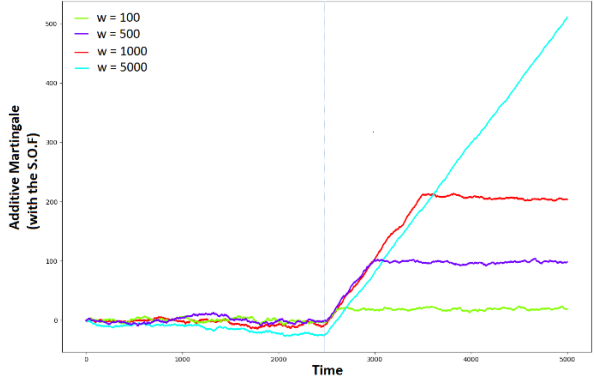

*Figure 2 (a): Additive martingale values over time with the shifted odd function as a betting function.*

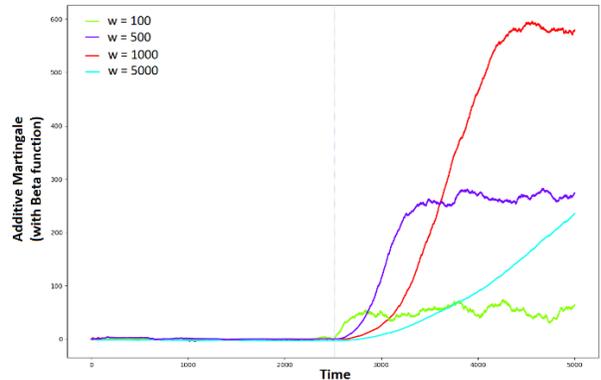

*Figure 2 (b): Additive martingale values over time with the estimated p.d.f. as a betting function.*

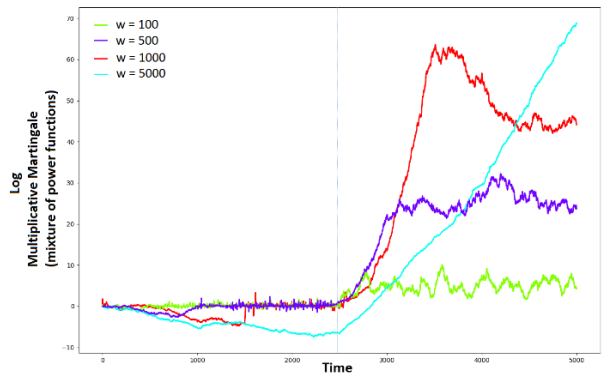

*Figure 2 (c): Multiplicative martingale values (in log scale) over time with the mixture of powers as a betting function.*

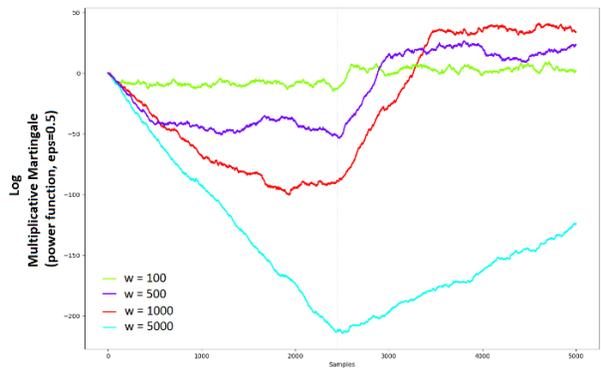

*Figure 2 (d): Multiplicative martingale values (in log scale) over time with the power function as a betting function.*

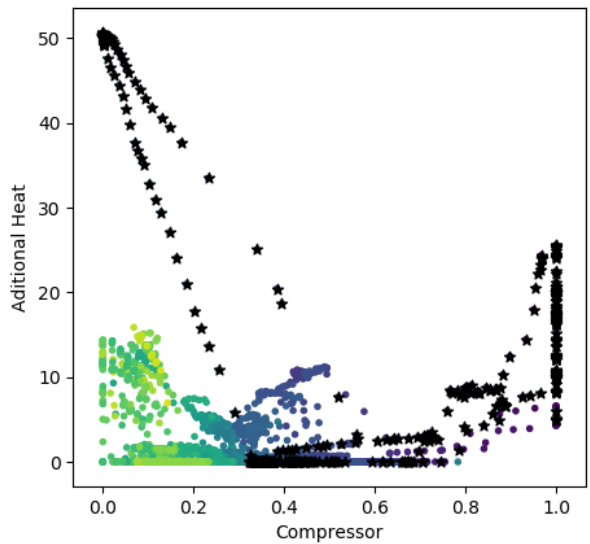
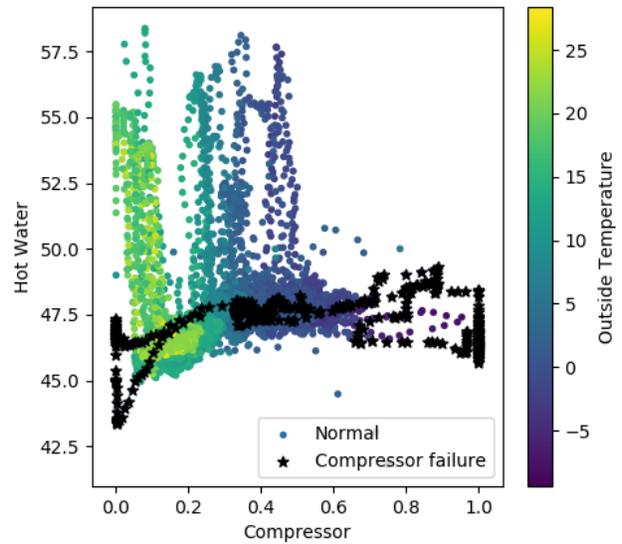

*Figure 3: Scatter plot of the heat-pump data over the whole period of time.*

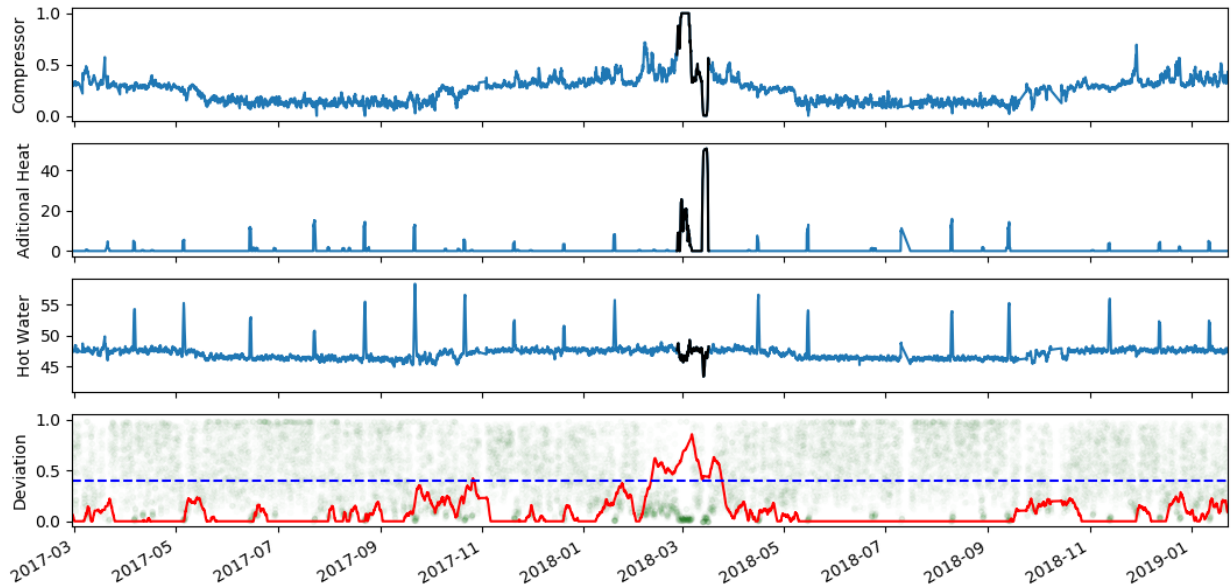

*Figure 4: Compressor failure detected on a real heat-pump using the proposed additive martingale method. The values of the sequence shown in red increased (above a threshold of 0.4) during the time where the compressor failure took place (shown in black on the time series plots).*

## 7. CONCLUSION AND FUTURE WORK

In this paper, a new additive martingale framework was proposed for detecting changes in streaming data based on the violation of exchangeability assumption. The proposed framework gives an alternative way for change-point detection, with advantages over existing ones. The proposed martingale sequence can be computed in a fully incremental way, it does not require any knowledge about the distribution of the data, and it enables theoretical guarantees on rejecting/accepting the null hypothesis "no change-point appears in the stream", using the Hoffding-Azuma inequality.

Nonetheless, the proposed method is not perfect and there are still many questions left open for further research: 1) A limitation of the current method is that it is not is not able to use an adaptive sliding window size. Indeed, as discussed in previous section, the change of martingale sequence in Figure 2 (b) is slow in the initial time steps when the size of sliding window is large. However, as can be observed from the same figure, when the sliding window size is smaller, the curve changes more rapidly in response to the deviation. This observation leads us to consider designing adaptive strategies, for example by using a smaller-size sliding window in the initial steps to increase the response speed to change, and gradually increasing the window size to gain more accurate information about the changed distribution in order to have a larger one-step increment; 2) There exist improvements over the basic Hoffding-Azuma inequality, for example some are presented in chapter 2 of **[Maxim 2015]**. It will be interesting to see whether these more advanced concentration inequalities can give tighter bounds than the one in eq. (21) for a given significance level.

## APPENDIX

The update rule for the sample mean is given as

$$\overline{p_n^w} = \overline{p_{n-1}^w} + \frac{1}{W}(p_n - p_{n-W}).$$

For sample variance, we have

$$M_n^w = \sum_{i=n-W+1}^{n} p_i^2 - W\left(\overline{p_n^w}\right)^2,$$

$$M_{n-1}^w = \sum_{i=n-W}^{n-1} p_i^2 - W\left(\overline{p_{n-1}^w}\right)^2,$$

which gives that

$$W\left(M_n^w - M_{n-1}^w - (p_n^2 - p_{n-M}^2)\right) = -\left(2\sum_{i=n-M+1}^{n-1} p_i\right)(p_n - p_{n-W}) - (p_n^2 - p_{n-W}^2)$$
$$= -W\left(\overline{p_n^w} + \overline{p_{n-1}^w}\right)(p_n - p_{n-W}),$$

with which we can conclude that

$$M_n^w = M_{n-1}^w + \left(p_n + p_{n-M} - \overline{p_n^w} - \overline{p_{n-1}^w}\right)(p_n - p_{n-W}).$$